\documentclass{amsart}
\usepackage{graphicx}

\newtheorem{thm}{Theorem}
\newtheorem{lemma}[thm]{Lemma}

\newcommand{\Nt}{\noindent {\bf Notation: }}
\newcommand{\Pf}{\noindent {\bf Proof: }}
\newcommand{\Rmk}{\noindent {\bf Remark: }}

\begin{document}

\title{Intrinsic Knotting and Linking of Almost Complete Graphs}
%\author{Thomas Mattman, Joel Pyzer,\\ Matt Rodrigues, Sam Williams}
\author[Campbell, Mattman, Ottman, Pyzer, Rodrigues, and Williams]{Jesse Campbell, Thomas W. Mattman, Ryan Ottman, 
Joel Pyzer, Matt Rodrigues, and Sam Williams}
%\maketitle

%\author{Thomas W. Mattman, Ryan Ottman, and Matt Rodrigues}        
\address{Department of Mathematics and Statistics,
         California State University, Chico,
         Chico CA 95929-0525, USA}
\email{Jesse.Austin.Campbell@gmail.com}
\email{TMattman@CSUChico.edu}
\address{Department of Mathematics,
	 South Hall Room 6607,
	 University of California,
Santa Barbara, CA 93106, USA}
\email{rottman@math.ucsb.edu}
\address{Middleton High School, 20932 Big Canyon Road, Middletown, CA 95461, USA}
\email{JDPyzer@juno.com}
\address{Department of Mathematics, UC Berkeley,
	 970 Evans Hall \#3840,
	 Berkeley, CA 94720-3840, USA}
\email{matthewrodrigues@msn.com}
\address{Department of Mathematics, UNC Asheville, CPO \#2350, One University Heights, Asheville, NC
28804-8511, USA}
\email{SBWillia@bulldog.unca.edu}
\thanks{
This work was completed in part through REU's at CSU, Chico during the
summers of 2003 to 2005 supported by NSF REU award \#0354174 and
the MAA's Strengthening Underrepresented Minority Mathematics Achievement
Program, with funding from the NSA, the
NSF, and the Moody's Foundation. The CSU, Chico Research Foundation provided additional support.}

\begin{abstract} 
We introduce new sufficient conditions for intrinsic knotting and
linking. A graph on $n$ vertices with at least $4n-9$ edges is
intrinsically linked. A graph on $n$ vertices with at least $5n-14$ edges
is intrinsically knotted. We also classify graphs that are $0$, $1$, or
$2$ edges short of being complete partite graphs with respect to intrinsic
linking and intrinsic knotting. In addition, we classify intrinsic
knotting of graphs on $8$ vertices. 
\end{abstract}

\maketitle

\section{Introduction}
We say that a graph is intrinsically knotted (respectively, linked)
if every tame embedding of the graph in $\mathbb{R}^3$ contains a non-trivially
knotted cycle (respectively, pair of non-trivially linked cycles).
Robertson, Seymour, and Thomas~\cite{RST}
demonstrated that intrinsic linking is determined by the seven
Petersen graphs. A graph is
intrinsically linked if and only if it  contains one of the  seven as a
minor. Since knotless embedding is preserved under edge contraction~\cite{NT}, work
of Robertson and Seymour~\cite{RS} shows that a similar, finite list
of graphs exists for the intrinsic knotting property. However,
determining this list remains difficult.
It is known~\cite{CG,F1,KS,MRS} that $K_7$ and
$K_{3,3,1,1}$ along with any graph obtained from these two by triangle-Y
exchanges is minor minimal with respect to intrinsic knotting.
Recently Foisy~\cite{F2,F3} has shown the existence of several additional
minor minimal intrinsically knotted graphs.

We will say that a graph is {\bf $k$-deficient} 
if it is a complete or complete partite graph
with $k$ edges removed. In the current
article we classify intrinsic knotting on
$1$- and $2$-deficient graphs as well as graphs on $8$ or fewer
vertices. As each of the intrinsically knotted graphs in
these families contains as a minor $K_7$, $K_{3,3,1,1}$, or a graph 
obtained from these by triangle-Y exchanges, we introduce no new examples
of minor minimal intrinsically knotted graphs.

On the other hand, we do introduce a new sufficient condition for 
intrinsic knotting:
\begin{thm}
\label{thmIK}
 A graph on $n$ vertices (where $n \geq 7$) 
with at least $5n-14$ edges is intrinsically knotted. 
\end{thm}
Moreover, we prove a similar condition for intrinsic linking.
\begin{thm} 
\label{thmIL}%
A graph on $n$ vertices (where $n \geq 6$) with at least
$4n-9$ edges is intrinsically linked.
\end{thm}
\Rmk Theorem~\ref{thmIL} was conjectured by Sachs~\cite{S} in 1984.

\bigskip

Theorems~\ref{thmIK} and \ref{thmIL} follow immediately from the
following 
\begin{thm}
\label{thmK7}
 A graph on $n$ vertices (where $n \geq 7$) 
with at least $5n-14$ edges has $K_7$ as a minor. 
\end{thm}
\begin{thm} 
\label{thmK6}%
A graph on $n$ vertices (where $n \geq 6$) with at least
$4n-9$ edges has $K_6$ as a minor.
\end{thm}
Indeed, $K_7$ is intrinsically knotted~\cite{CG} and if a graph $G$
has an intrinsically knotted minor, then $G$ is intrinsically knotted~\cite{NT}.
Similarly, as $K_6$ is intrinsically linked~\cite{CG, S}, any graph
containing it as a minor is likewise intrinsically linked.

Motivated by Sachs's Conjecture, we found a proof of
Theorem~\ref{thmK6}. Recently, we learned that Mader~\cite{M} had 
proved both Theorems~\ref{thmK7} and \ref{thmK6} in 1968! As our
proof of Theorem~\ref{thmK6} is quite different from Mader's, we
include it as Section 2 of this paper. 

In Sections 3, 4, and 5 we
investigate $0$-, $1$-, and $2$-deficient graphs respectively. In each
case we classify the graphs with respect to intrinsic knotting and
linking. (The classification of intrinsic knotting of $0$-deficient
graphs is due to Blain, Bowlin, Fleming, Foisy, Hendricks, and
LaCombe~\cite{Fl}.) In Section 6 we classify intrinsically knotted
graphs on $8$ vertices. (This classification  was recently
independently presented in \cite{Fl}.)

\section{Proof of Theorem~\ref{thmK6}}

\setcounter{thm}{3}

In this section we prove:

\begin{thm} 
A graph on $n$ vertices (where $n \geq 6$) with at least
$4n-9$ edges has $K_6$ as a minor.
\end{thm}

We begin with some definitions.
A {\bf graph} $G$ is a finite vertex set $V(G)$ and edge family $E(G)$,
and we mean $G$ to be simple, connected, and undirected.
A graph is labeled
$K_{n}$, and called the {\bf complete graph} on $n$ vertices, if every
vertex in the graph is adjacent to every other vertex.
A graph $H$ is a {\bf minor} of a
graph $G$ if $H$ results from a finite series of edge contractions,
vertex deletions, and edge deletions on $G$.
Two adjacent vertices
$x$ and $y$ have a {\bf common neighbor} $z$ in $G$ if
$xz,yz\in E(G)$.
%Let $CN(xy)$ be the set of common
%neighbors between two adjacent vertices $x$ and $y$, i.e.,
%$CN(xy)=\{ v\in V(G)|vx,vy\in E(G)\}$.

\bigskip

\Pf (of Theorem~\ref{thmK6}) 
We proceed using induction on $n$, the number of vertices of our
graph $G$. If $n=6$, then $G$ has at least $4(6)-9=15$ edges. Since
$K_6$ is the only graph on $6$ vertices with $15$ edges, $G = K_6$.

We now assume that any graph on $6,7,...n-1$ vertices  (with at least 
$4 | V(G) |-9$ edges) has a
$K_6$ minor. Let $G$ be a graph on $n$ vertices that has $4n-9$ edges.
(Note that if
$\left| E(G)\right| >4n-9$, then $G$ has a minor on $4n-9$ edges which
can be obtained through edge deletions.)

If $G$ has two adjacent vertices with less than four
common neighbors, then contracting their edge results in
a graph on $n-1$ vertices with at least $4n-9 - 4 = 4(n-1) -9$ edges.
By induction, $G$ has a $K_6$ minor. 

Therefore, we may assume
that each pair of adjacent vertices in $G$ has at least four common
neighbors.
This implies that every vertex in $G$ has degree at least five. 
Suppose every vertex in $G$ has degree eight or more. Since the number
of edges is half the degree sum, this would mean $G$ has at least 
$4n$ edges contradicting our assumption of $4n-9$ edges. The remainder of
the argument will focus on a vertex $a$ of least degree. That is, $a\in
V(G)$ and
$\forall v\in V(G)$, $\mbox{deg}(a)\leq \mbox{deg}(v)$. We are left with three cases,
$\mbox{deg}(a) = 5$, $6$, and $7$.

Suppose  $\mbox{deg}(a) = 5$. Consider the subgraph $H$ on
the neighbors of $a$.  Each of the neighbors of $a$ must share at least
four common neighbors with it and, therefore, has degree $4$ in $H$. 
But then $H$ is a simple graph with
five vertices and $\frac{5(4)}{2} =10$ edges, i.e., $H = K_5$. Thus,
because $a$ is connected to every vertex in $H$ , the graph on
$a$ and its neighbors is $K_6$ , implying $G$ has a $K_6$ minor.

\begin{figure}[h]
\begin{center}
\includegraphics[scale=.3]{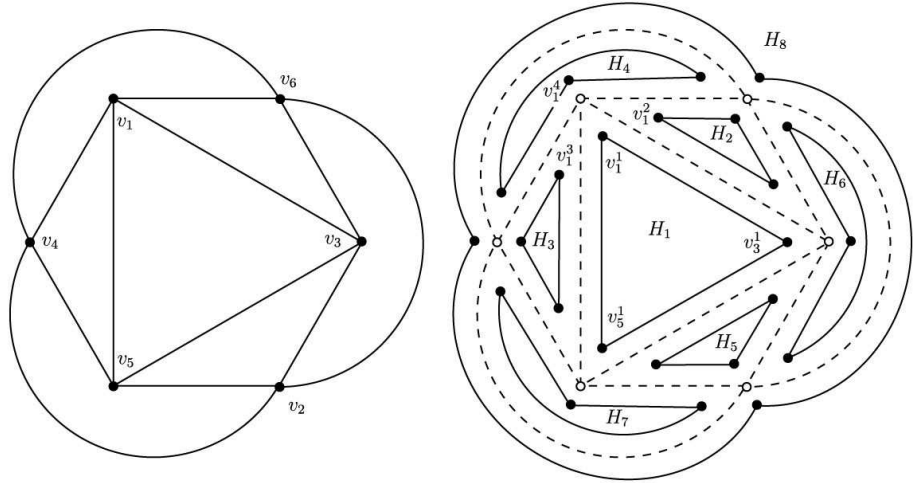}
\caption{\label{figdeg6}%
A vertex of degree $6$}
\end{center}
\end{figure}

Now, suppose $\mbox{deg}(a)=6$. Consider the graph $H$
on the neighbors of $a$: $v_1$, $v_2$, $v_3$, $v_4$, $v_5$, and $v_6$.
Since every vertex $v_{i}$ must share at least four common neighbors with
$a$, $\forall v_{i}\in V(H), \mbox{deg}_{H}(v_{i})\geq 4$. It follows that
$H$ has at least $\frac{6(4)}{2} = 12$ edges.
However, we know that if $\left|E(H)\right| \geq 3 |V(H)|-5=3(6)-5=13$, then
$H$ contains a $K_5$ minor \cite[Corollary 1.13]{B}, implying that $G$
contains a $K_6$ minor. So, we may assume that $H$ has exactly $12$ edges
and each vertex has degree $4$. There is only one such graph: a degree
4-regular triangulation on six vertices (see left of
Figure~\ref{figdeg6}).
Notice that $H$ is $K_6$ less the  edges $v_1 v_2$, $v_3 v_4$, and $v_5
v_6$. If we add back any of those three edges $H$ will have at least 13
edges, hence a $K_5$ minor and then $G$ will have a $K_6$ minor. 

Let $G'$ be the subgraph formed from $G$ by removing $a$ and its six
edges. Then $G'$ has $n-1$ vertices and $4(n) -9 - 6 = 4n-15$ edges.
Let $G' \setminus H$ denote the graph obtained by removing the 
edges of $H$. We can assume that there is no path in $G' \setminus H$
from $v_1$ to $v_2$ , for otherwise, we can contract
along the edges of the path to form a minor in $G$ of the form
$a + H'$ with $H'= H \cup v_1 v_2$ having 13 edges. That is, $G$ would
have a $K_6$ minor. Similarly,
in $G' \setminus H$ we can assume there is no path between $v_3$ and
$v_4$, nor any path linking $v_5$ and $v_6$. This means that each vertex
in $V(G') \setminus V(H)$ is connected to at most  three of the vertices
in $\{ v_1, \ldots v_6 \}$ by edges in $G' \setminus H$. 

This will allow us to divide $G'$ into eight subgraphs that overlap only in the subgraph $H$. 
To illustrate this, we will construct a new graph $G''$ consisting of the disjoint union of the
eight subgraphs. 
Let's first use $H$  to create eight disjoint triangles $H_1, \ldots, H_8$, as at right in Figure 1.
Note that this entails doubling each edge of $H$ and making four copies of each vertex. The vertices of
$G''$ will be those of $V(G') \setminus V(H)$ along with the vertices of $H_1, \ldots H_8$.
So, we've added three extra copies of each of the six vertices in $H$ and, therefore, $G''$ will have 
$n-1 + 3(6) = n+17$ vertices. The edges of $G''$ will be those of $G' \setminus H$ along with two
copies of each  edge of $H$. Thus, $G''$ will have $4n - 15 + 12 = 4(n+17) - 71$ edges.

Now, the vertices of $H_1, \ldots, H_8$ are connected as eight $K_3$'s and
the vertices in $V(G'') \setminus ( V(H_1) \cup \ldots  \cup V(H_8))$ are connected to one another as they were in
$G' \setminus H$. To complete the construction of $G''$, we need only describe how the edges between vertices 
of $H$ and vertices in 
$V(G') \setminus V(H)$ are placed in $G''$. For example, let's look at the edges incident on $v_1$.
Suppose $w \in V(G') \setminus V(H)$  is a neighbor of $v_1$. Then $w$ can be connected in $G' \setminus H$
to at most two other vertices of $H$. For example, if $w$ is also connected to $v_3$ and $v_5$, then in
$G''$, we will  connect  $w$ to the vertex $v_1^1$ in $H_1$, because $H_1$ is the triangle
corresponding to the  triangle $v_1 v_3 v_5$ in $H$. If $w$ is connected to at most one other vertex in
$H$, say $v_3$,  then in $G''$,  we have the choice of connecting $w$ to $H_1$ or $H_2$ as
those are the two triangles containing copies of the vertices $v_1$ and $v_3$.
We will adopt the convention of connecting it to the $H_i$ of lower index, in
 this case $H_1$. That is, in $G''$ there will be an edge between $w$ and the vertex $v_1^1$ in
$H_1$. Similarly, if $w$ is connected to no other vertex in $H$, we have four choices, $H_1, H_2, H_3,
H_4$, of where to  connect  $w$ in $G''$. Again, we will  opt for
$H_1$, the choice of  lowest index. In a similar way, each  neighbor $w$ of $v_2, \ldots, v_6$  is
attached to  one of the  graphs $H_i$. This placement is unambiguous if $w$  is connected to three
vertices of $H$, it is made  among a choice of two $H_i$ if $w$ is connected to only two vertices,  and
it is made from a choice of  four $H_i$ if $w$ is connected to only one vertex in $H$. When there is a
choice, we will use the $H_i$ of  least index $i$ among the given options. 

This completes the construction of $G''$. Let $G_i$ be the connected component of $H_i$ in $G''$. 
Since $G$ was connected,  $G'' \subset \cup_{i=1}^8 G_i$. On the other hand, since we are assuming 
there is no path in $G' \setminus H$ connecting $v_1$ and $v_2$, $v_3$ and $v_4$, or $v_5$ and $v_6$,
the $G_i$ must be disjoint. 
Now, $G''$ has $n+17$ vertices and $4(n+17) - 71$ edges. By the pigeon hole principle, this means that
one of the subgraphs $G_i$ has $n_i$ vertices and at least $4n_i - \lfloor 71/8 \rfloor = 
4n_i - 8$. edges. This means $G_i$ is not simply a triangle, so $n_i \geq 7$ (recall each vertex has degree
at least six).  Thus, $G_i$ is a proper subgraph of $G$ that has a $K_6$ minor. 

\begin{figure}[ht]
\begin{center}
\includegraphics[scale=.3]{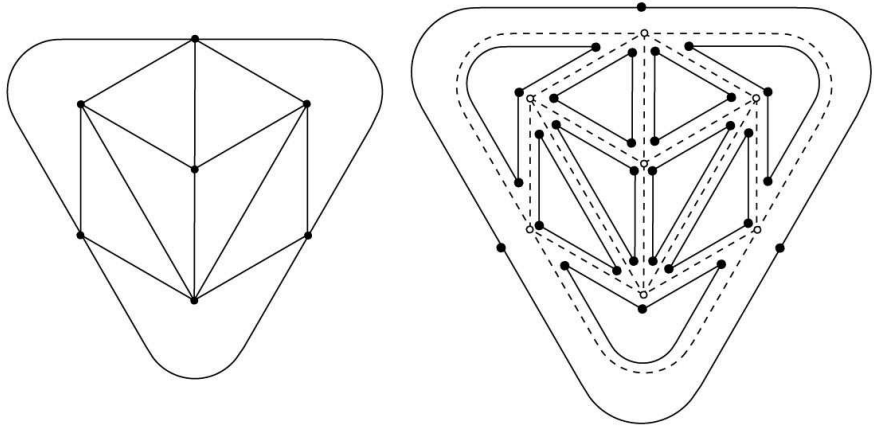}
\caption{\label{figdeg7}%
A vertex of degree $7$.}
\end{center}
\end{figure}

Finally, suppose $\mbox{deg}(a) = 7$. Let $H$ be the induced graph on the neighbors of $a$. Again, 
each vertex must have degree $4$ or more and there are fewer than $3(7) - 5 = 16$ edges. 
There are several such graphs, but only one has no $K_5$ minor.
This time it is a triangulation on seven
vertices with two vertices of degree five and five vertices of degree
four, and with 10 triangles (see Figure~\ref{figdeg7}). In this case $G''$ has $m=n+22$ vertices and
$4m-89$ edges. There must be a subgraph of $G''$ that has at least
$4n_{i}-\left\lfloor \frac{89}{10}\right\rfloor \geq 4n_{i}-8$ edges. So,
$G$ has a subgraph that has a $K_6$ minor, implying $G$ has a $K_6$
minor. \qed

\section{Complete and Complete Partite Graphs}
In this section we present a classification of intrinsic linking of $0$-deficient graphs, i.e., complete
and complete partite graphs. We also
include the classification of intrinsic knotting found in \cite{Fl} for the reader's convenience.
We begin with some lemmas.

If we know that a graph $G$ contains a linked (knotted) minor, then $G$ 
must also be linked (knotted)~\cite{NT}. Conversely, if we can realize $G$
as a minor of an unlinked (unknotted) graph, we deduce that
$G$ must also be unlinked (unknotted). A useful lemma in this 
regard shows how we can combine parts of a $k$-deficient graph. Recall that 
$K_{n_1, n_2, \ldots, n_p}$ denotes the complete $p$ partite graph
with $n_i$ vertices in the $i$th part. Let $K_{n_1,n_2, \ldots, n_p} - k$ denote a graph obtained by removing $k$ edges from $K_{n_1,n_2, \ldots, n_p}$.

\begin{lemma}\label{combparts}
$K_{n_1+n_2, n_3, \ldots, n_p}$ is a minor of $K_{n_1,n_2, \ldots, n_p}$.
Similarly, every graph of the form
$K_{n_1,n_2, \ldots, n_p} - k$ has a minor of the form
$K_{n_1+n_2, n_3, \ldots, n_p} - k$.
\end{lemma}

\Pf
Combining the $n_1$  and $n_2$ parts only involves removing edges between
vertices in the $n_1$ part and the $n_2$ part. Recall that with complete
partite graphs, the ordering of the subscripts is not important, so this
lemma implies that we can combine any 2 parts  to get a minor of
the original graph. 

Now, suppose we have a complete partite graph with $k$ edges removed
and we combine two parts. The result would be a complete partite
graph with one fewer part and at most $k$ edges removed.  

Furthermore, if, in $K_{n_1,n_2, \ldots, n_k}-k$, there are $m$ edges missing between parts $n_1$ and $n_2$, we can see that 
$K_{n_1,n_2, \ldots, n_k}-k$ has a minor of the form
$K_{n_1+n_2, n_3, \ldots, n_k}-(k-m)$. 
\qed

\medskip

The other main technique we will use in our classifications (aside from
Theorems~\ref{thmIK} and \ref{thmIL} and Lemma~\ref{combparts})
is based on a lemma for intrinsic linking due
to Sachs~\cite{S} and an analogous result for intrinsic knotting
due to \cite{Fl,OT}. Let $G + H$ denote the join of graphs $G$ and
$H$, i.e., the graph obtained by taking the union of $G$ and $H$ and
adding an edge between each vertex
of $G$ and each vertex of $H$.

\begin{lemma}[\cite{S}]\label{sachs}

The graph $G + K_1$ 
is intrinsically linked if and only if G is non-planar

\end{lemma}

%\begin{cor} 
%If a vertex is deleted from a graph $H$ and the result is a
%planar graph, then $H$ is not intrinsically linked. 
%\end{cor}

\begin{lemma}[\cite{Fl,OT}]\label{flemming}

The graph $G + K_2$ is intrinsically knotted if and only if
$G$ is non-planar

\end{lemma}

%\begin{cor}
%If two vertices are deleted from a graph $H$ and the result is a
%planar graph, then $H$ is not intrinsically knotted.
%\end{cor} 

Surveying the methods that we will use in the rest of the paper, 
we see that we will not be introducing any new minor minimal graphs 
for intrinsically knotting. Indeed,
if we use Lemma~\ref{flemming} to argue that a graph 
$H = G+K_2$ is
intrinsically knotted, then, $G$ must contain
$K_5$ or $K_{3,3}$ as a minor. It follows that $H$ contains
$K_7$ or $K_{3,3,1,1}$ as a minor. 
Similarly, Theorem~\ref{thmIK} is based on showing that the
graph has a $K_7$ minor. As our only other technique (demonstrating an 
intrinsically knotted minor)
relies on the graphs derived from $K_7$ and $K_{3,3,1,1}$ by 
triangle-Y exchanges,  our methods cannot introduce new 
minor minimal examples. In other words,
every intrinsically knotted graph that is $0$-, $1$-, or $2$-deficient
or a graph on $8$ vertices necessarily has as a minor $K_7$, 
$K_{3,3,1,1}$ or a graph obtained from one of these by
triangle-Y exchanges.

\begin{thm}

The complete $k$-partite graphs are classified with respect to intrinsic 
linking according to Table \ref{ilcomp}.
\end{thm}

\begin{table}[h]
\begin{center}
\begin{tabular}{|l||c|c|c|c|c|c|} \hline k          & 1 &  2  & 3     &
4       & 5         & $\geq$ 6 \\ \hline  linked     & 6 & 4,4 & 3,3,1 &
2,2,2,2 & 2,2,1,1,1 & All \\ 
           &   &     & 4,2,2 & 3,2,1,1 & 3,1,1,1,1 & \\ \hline not 
linked     & 5 & $n$,3 & 3,2,2 & 2,2,2,1 & 2,1,1,1,1 & None \\
           &   &   & $n$,2,1 & $n$,1,1,1 &           & \\ \hline
\end{tabular}
\end{center}
\caption{Intrinsic linking of complete $k$-partite graphs.}\label{ilcomp}
\end{table}

\Rmk The $n$ in the notation $K_{n,3}$ indicates that the property holds for any
number of vertices in the first part, i.e., none of the graphs $K_{1,3}$, $K_{2,3}$,
$K_{3,3}$, $\ldots$ are intrinsically linked.
For each $k$,  the table includes minimal examples of
intrinsically linked complete $k$-partite graphs and maximal graphs which
are not intrinsically linked.
For example, any complete $2$-partite graph which contains
$K_{4,4}$ as a minor is intrinsically linked. On the other hand, any
complete $2$-partite graph which is a minor of a $K_{n,3}$ 
is not linked. Thus,
$K_{l,m}$ is intrinsically linked if and only if $l \geq 4$ and $m \geq 4$.

\Pf

\noindent%
Let us demonstrate that the graphs labeled ``linked" in Table
\ref{ilcomp} are in fact intrinsically linked. 

Sachs~\cite{S}, and, independently, Conway and Gordon~\cite{CG} proved that
$K_6$ is linked. Any $k$-partite graph with $k \geq 6$ contains  $K_6 =
K_{1,1,1,1,1,1}$ as a minor and is therefore linked. For the remaining
$k$, we appeal to  work of Robertson, Seymour, and Thomas
~\cite{RST} who showed that a graph   is intrinsically linked if and only if
it  contains as a minor one of the seven graphs in the Petersen family. In
particular, $K_{4,4}$ with one edge removed and $K_{3,3,1}$ are both in
this family. By combining parts, we see that, for $2 \leq k \leq 5$, one
of these two is a minor of each of the ``linked" graphs in Table 1. 

For each of the ``not linked" examples in Table 1 which involve a part
with a single vertex, the corresponding graph obtained by removing that
vertex is planar. Therefore, by Lemma \ref{flemming}, these graphs are
not intrinsically linked. 

Since a cycle requires at least three  vertices, $K_5$ has no disjoint
pair of cycles and is therefore not linked. The remaining ``not linked" 
graphs in Table 1, $K_{n,3}$ and $K_{3,2,2}$, are, respectively, minors
of the unlinked graphs $K_{n,2,1}$ and
$K_{2,2,2,1}$.
\qed

\medskip

For the reader's convenience, we reproduce the 
classification of knotted partite graphs due to Blain, Bowlin, Fleming, Foisy, Hendricks,
and LaCombe~\cite{Fl} as Table \ref{ikcomp} below.

\begin{table}[h]
\begin{center}
\begin{tabular}{|l||c|c|c|c|c|c|c|} \hline k          & 1 &  2  & 3     &
4       & 5         & 6 & $\geq 7$ \\
\hline   knotted    & 7 & 5,5 & 3,3,3 & 3,2,2,2 & 2,2,2,2,1 & 2,2,1,1,1,1
& All \\ 
           &   &     & 4,3,2 & 4,2,2,1 & 3,2,2,1,1 & 3,1,1,1,1,1 &\\
           &   &     & 4,4,1 & 3,3,1,1 & 3,2,1,1,1 &  & \\ \hline

not  knotted    & 6 & 4,4 & 3,3,2 & 2,2,2,2 & 2,2,2,1,1 & 2,1,1,1,1,1 &
None \\
           &   &   & $n$,2,2 & 4,2,1,1 & 2,2,1,1,1 & &\\
           &   &   & $n$,3,1 & 3,2,2,1 & $n$,1,1,1,1 & & \\ 
           &   &   &         & $n$,2,1,1 &         & & \\ \hline
\end{tabular}
\end{center}
\caption{Intrinsic knotting of $k$-partite graphs.}\label{ikcomp}
\end{table}

\section{$1$-deficient graphs}

In this section we classify $1$-deficient graphs with respect to 
intrinsic linking and knotting.  

\Nt Often, we will have to talk about a particular vertex or part. We
will refer to parts alphabetically with capital letters and vertices of
those parts with lower case letters. For example, in
$K_{4,3,1}$, we will call the part with 4 vertices part A, the part with
3 vertices part B, and the part with 1 vertex part C. An edge between
parts A and C would be labeled (a,c). 

\subsection{Intrinsic linking}

\begin{thm}

The 1-deficient graphs are classified with respect to intrinsic  linking
according to Table \ref{il1def}.
\end{thm}

\begin{table}[ht]
\begin{center}
\begin{tabular}{|l||c|c|c|c|c|c|c|} \hline k          & 1   &  2      &
3       & 4             & 5               & 6             & $\geq$ 7 \\
\hline  linked     & 7-e & 4,4-e   & 4,3,1-e & 2,2,2,2-e     &
2,2,1,1,1-(b,c) & 2,1,1,1,1,1-e & All \\ 
           &     &         & 3,3,2-e & 3,2,1,1-(b,c) & 3,1,1,1,1-(b,c)
&              & \\  
           &     &         & 4,2,2-e & 4,2,1,1-e     & 4,1,1,1,1-e    
&              & \\  
           &     &         &         & 3,3,1,1-e     & 3,2,1,1,1-e    
&              & \\ 
           &     &         &         & 3,2,2,1-e     & 2,2,2,1,1-e    
&              & \\ \hline

not        & 6-e & $n$,3-e & 3,2,2-e & 2,2,2,1-e     &
2,2,1,1,1-(a,b) & 1,1,1,1,1,1-e & None \\
linked     &     &       & $n$,2,1-e & $n$,1,1,1-e   & 2,2,1,1,1-(c,d)
&              & \\  
           &     &       & 3,3,1-e   & 3,2,1,1-(a,b) & 3,1,1,1,1-(a,b)
&              & \\  
           &     &       &           & 3,2,1,1-(a,c) & 2,1,1,1,1-e    
&              & \\  
           &     &       &           & 3,2,1,1-(c,d) &                
&              & \\ \hline

\end{tabular}
\end{center}
\caption{Intrinsic linking of 1-deficient graphs.}\label{il1def}
\end{table}

\Rmk

Note that most of these graphs have different types of edges. 
For some graphs, removal of one edge (e.g., (a,b)) will result in a
non-linked graph while the removal of a different edge (e.g., (b,c)) will
result in a linked graph. However, in many cases, removal of any edge
will result in the same categorization. In such cases we simply write
``-e''. For example, no matter which edge we remove from $K_{3,2,2}$, we will
obtain a graph that is not intrinsically linked. 

\Pf

%\noindent%
%\underline{Linked:}

Let us demonstrate that the graphs labeled ``linked" in Table
\ref{il1def} are in fact intrinsically linked. 

By Theorem~\ref{thmIL}, $K_7$-e, $K_{2,2,2,2}$-e, $K_{3,2,1,1,1}$-e,
$K_{2,2,2,1,1}$-e, and $K_{2,1,1,1,1,1}$-e are intrinsically linked. Any $1$-deficient graph with
seven or more parts contains $K_7$-e as a minor and is therefore intrinsically linked.

Recall that $K_{4,4}$-e is a Petersen graph. By Lemma~\ref{combparts}, it's a minor
of $K_{4,3,1}$-e,
$K_{4,2,2}$-e, $K_{4,2,1,1}$-e, $K_{3,3,1,1}$-e,
$K_{3,2,2,1}$-e, and $K_{4,1,1,1,1}$-e.

$K_{3,3,2}$-e, and, by Lemma~\ref{combparts}, 
 $K_{3,2,1,1}$-(b,c), $K_{2,2,1,1,1}$-(b,c), and
$K_{3,1,1,1,1}$-(b,c) all contain the Peterson graph $K_{3,3,1}$ as a minor and are
therefore intrinsically linked. 

\begin{table}[hb]
\begin{center}
\begin{tabular}{|l||c|c|c|c|c|} \hline k          & 1   &  2    &  
3     &    4    &     5     \\
%&      6      &  7                &  $\geq 8$ \\
\hline   
knotted    & 8-e & 5,5-e & 3,3,3-e & 3,2,2,2-e & 2,2,2,2,1-e \\
%& 2,2,1,1,1,1-(b,c) & 2,1,1,1,1,1,1-e & All \\ 
           &     &       & 4,3,2-e & 4,2,2,1-e & 3,2,1,1,1-(b,c) \\
% & 3,1,1,1,1,1-(b,c) &        & \\ 
           &     &       & 4,4,1-e & 3,3,2,1-e & 4,2,1,1,1-e \\ 
% & 3,2,1,1,1,1-e &                   & \\ 
           &     &     &           & 4,3,1,1-e & 3,3,1,1,1-e \\
% & 2,2,2,1,1,1-e &                   & \\ 
           &     &     &           &           & 3,2,2,1,1-e \\
% & 4,1,1,1,1,1-e &                   & \\ 

\hline

%          &     &     &           &           &             & 
%&                   & \\ 

not  knotted    & 7-e & n,4-e & 3,3,2-e &  3,3,1,1-e & 3,2,1,1,1-(c,d) \\
% & 2,2,1,1,1,1-(a,b) & 1,1,1,1,1,1,1-e & None \\ 
           &     &     & $n$,2,2-e &  2,2,2,2-e & 3,2,1,1,1-(a,b) \\
% & 2,2,1,1,1,1-(c,d) &     & \\ 
           &     &     & $n$,3,1-e &  3,2,2,1-e & 3,2,1,1,1-(a,c) \\
% & 3,1,1,1,1,1-(a,b) &                   & \\ 
           &     &     &         &  $n$,2,1,1-e &  2,2,2,1,1-e  \\
% & 2,1,1,1,1,1-e &                   & \\ 
           &     &     &       &            &  $n$,1,1,1,1-e    \\
% &  &                   & \\ \hline

\hline
\multicolumn{6}{c}{ } \\
\hline

k          & \multicolumn{3}{c|}{6}                 & 7               & $\geq 8$ \\ \hline   
knotted    & \multicolumn{3}{c|}{2,2,1,1,1,1-(b,c)} & 2,1,1,1,1,1,1-e & All      \\ 
           & \multicolumn{3}{c|}{3,1,1,1,1,1-(b,c)} &                 &          \\ 
           & \multicolumn{3}{c|}{3,2,1,1,1,1-e}     &                 &          \\ 
           & \multicolumn{3}{c|}{2,2,2,1,1,1-e}     &                 &          \\ 
           & \multicolumn{3}{c|}{4,1,1,1,1,1-e}     &                 &          \\ \hline

%          &     &     &           &           &             & 
%&                   & \\ 

not  knotted & \multicolumn{3}{c|}{2,2,1,1,1,1-(a,b)} & 1,1,1,1,1,1,1-e & None \\ 
             & \multicolumn{3}{c|}{2,2,1,1,1,1-(c,d)} &                 &      \\ 
             & \multicolumn{3}{c|}{3,1,1,1,1,1-(a,b)} &                 &      \\ 
             & \multicolumn{3}{c|}{2,1,1,1,1,1-e}     &                 &      \\ \hline

\end{tabular}
\end{center}
\caption{Intrinsic knotting of 1-deficient graphs.}\label{ik1def}
\end{table}

\begin{table}[ht]

\begin{center}

\begin{tabular}{|l||c|c|c|c|} \hline

k          & 1    & 2      & 3                                & 4\\ \hline

linked     & 7-2e & 5,4-2e & 4,3,1-\{($b,c),e$\}                &
3,2,1,1-\{($b_1,c),(b_1,d$)\} \\ 

           &      &        & 4,3,1-\{($a_1,b_1),(a_1,b_2$)\}  &
3,2,1,1-\{$(b_1,c),(b_2,c)$\} \\ \cline{5-5}

           &      &        & 4,3,1-\{($a_1,b_1$),($a_1,c$)\}  &
4,2,1,1-\{$(b,c),e$\} \\ \cline{4-4}

           &      &        & 3,3,2-\{($a_1,b_1),(a_1,b_2$)\}  &
4,2,1,1-\{$(c,d),e$\} \\ 

           &      &        & 3,3,2-\{($a_1,c_1),(a_2,c_1$)\}  &
4,2,1,1-\{$(a_1,b_1),(a_1,b_2)$\} \\ 

           &      &        & 3,3,2-\{($a_1,c_1),(b_1,c_1$)\}  &
4,2,1,1-\{$(a_1,b_1),(a_1,c)$\} \\ 

           &      &        & 3,3,2-\{($a_1,b_1),(b_2,c_1$)\}  &
4,2,1,1-\{$(a_1,c),(a_1,d)$\} \\ \cline{5-5} 

           &      &        & 3,3,2-\{($a_1,c_1),(b_1,c_2$)\}  &
3,3,1,1-\{$(b,c),e$\} \\ 

           &      &        & 3,3,2-\{($a_1,c_1),(a_1,c_2$)\}   &
3,3,1,1-\{$(c,d),e$\} \\ \cline{4-4}

           &      &        & 4,2,2-\{$(b,c),e$\} &
3,3,1,1-\{$(a_1,b_1),(a_1,b_2)$\} \\ \cline{5-5}

           &      &        & 4,2,2-\{($a_1,b_1),(a_1,b_2$)\} &
3,2,2,1-\{$(b,c),e$\} \\ \cline{4-4}

           &      &        & 5,3,1-2e & 3,2,2,1-\{$(c,d),e$\} \\ 

           &      &        & 4,4,1-2e  & 3,2,2,1-\{$(a,d),e$\} \\ 

           &      &        & 4,3,2-2e   &
3,2,2,1-\{$(a_1,b_1),(a_1,b_2)$\} \\ 

           &      &        &  3,3,3-2e  &
3,2,2,1-\{$(a_1,b_1),(a_2,b_1)$\} \\ 

           &      &        &  5,2,2-2e &
3,2,2,1-\{$(a_1,b_1),(a_2,c_1)$\} \\ \cline{5-5}

           &      &        &   & 2,2,2,2-2e \\ 

           &      &        &                          & 5,2,1,1-2e \\ 

           &      &        &                          & 4,3,1,1-2e \\

           &      &        &                         & 4,2,2,1-2e \\

           &      &        &                         & 3,3,2,1-2e \\
\hline

%           &      &        &                          & 3,2,2,2-2e \\ 

%           &      &        &                                  & \\ 

not linked & 6-2e & 4,4-2e & 4,3,1-\{($a_1,b_1),(a_2,b_1$)\}  &
3,2,1,1-\{($a,b),e$\}\\ 

           &      & n,3-2e & 4,3,1-\{($a_1,b_1),(a_2,b_2$)\}  &
3,2,1,1-\{($a,c),e$\}\\ 

           &      &        & 4,3,1-\{($a_1,c),(a_2,c$)\}      &
3,2,1,1-\{($c,d),e$\}\\ 

           &      &        & 4,3,1-\{($a_1,b_1),(a_2,c$)\}    &
3,2,1,1-\{($b_1,c),(b_2,d)$\}\\ \cline{4-5} 

           &      &        & 3,3,2-\{($a_1,b_1),(a_2,b_2$)\}  &
4,2,1,1-\{$(a_1,b_1),(a_2,b_2)$\}\\ 

           &      &        & 3,3,2-\{($a_1,b_1),(b_1,c_1$)\}  &
4,2,1,1-\{$(a_1,c),(a_2,c)$\}\\ 

           &      &        & 3,3,2-\{($b_1,c_1),(b_2,c_2$)\}  &
4,2,1,1-\{$(a_1,b_1),(a_2,b_1$)\} \\ \cline{4-4}

           &      &        & 4,2,2-\{($a_1,b_1),(a_2,b_2$)\}  &
4,2,1,1-\{$(a_1,b_1),(a_2,c)$\}\\ 

           &      &        & 4,2,2-\{($a_1,b_1),(a_2,b_1$)\}  &
4,2,1,1-\{$(a_1,c),(a_2,d)$\}\\ \cline{5-5}

           &      &        & 4,2,2-\{($a_1,b_1),(a_2,c_1$)\}  &
3,3,1,1-\{$(a_1,b_1),(a_2,b_2$)\}\\ \cline{5-5}

           &      &        & 4,2,2-\{($a_1,b_1),(a_1,c_1$)\}  &
3,2,2,1-\{$(a_1,b_1),(a_2,b_2$)\}\\ \cline{4-4}

           &      &        & 3,2,2-2e                         &
3,2,2,1-\{$(a_1,b_1),(a_1,c_1)$\}\\ \cline{5-5}

           &      &        & n,2,1-2e                         &
2,2,2,1-2e\\ 

           &      &        & 3,3,1-2e                         &
n,1,1,1-2e\\ \hline

\end{tabular}

\end{center}

\caption{Intrinsic linking of 2-deficient graphs}\label{il2def1}

\end{table}

\begin{table}[ht]

\begin{center}

\begin{tabular}{|l||c|c|c|} \hline

k           & 5                                   &
6                                 & $\geq$7 \\ \hline

linked      & 2,2,1,1,1-\{$(b_1,c),(b_2,c)$\}     &
2,1,1,1,1,1-\{$(a_1,b),(a_1,c)$\}  & All \\ 

            & 2,2,1,1,1-\{$(a_1,d),(b_1,c)$\}     &
2,1,1,1,1,1-\{$(a_1,b),(a_2,b)$\}  &  \\ 

            & 2,2,1,1,1-\{$(b_1,c),(b_1,d)$\}     &
2,1,1,1,1,1-\{$(a_1,b),(c,d)$\} &  \\ \cline{2-2}

            & 3,1,1,1,1-\{$(b,c),(c,d)\}$         &
2,1,1,1,1,1-\{$(b,c),(c,d)$\} &  \\ \cline{2-3}

            & 4,1,1,1,1-\{$(b,c),e$\}             & 3,1,1,1,1,1-2e &  \\

            & 4,1,1,1,1-\{$(a_1,b),(a_1,c)$\}     & 2,2,1,1,1,1-2e &  \\
\cline{2-2}

            & 3,2,1,1,1-\{$(a,c),e$\}             &  &  \\

            & 3,2,1,1,1-\{$(b,c),e$\} &  &  \\

            & 3,2,1,1,1-\{$(c,d),e$\} &  &  \\

            & 3,2,1,1,1-\{$(a_1,b_1),(a_1,b_2)$\} &  &  \\

            & 3,2,1,1,1-\{$(a_1,b_1),(a_2,b_1)$\} &  &  \\ \cline{2-2}

            & 2,2,2,1,1-2e                        &  &  \\

            & 5,1,1,1,1-2e                        &  &  \\

            & 4,2,1,1,1-2e                        &  &  \\

            & 3,3,1,1,1-2e                        &  &  \\ \hline

 %           & 3,2,2,1,1-2e                        &  &  \\

%            &  &  &  \\

not linked  & 2,2,1,1,1-\{$(a,b),e$\}             &
2,1,1,1,1,1-\{$(a_1,b),(a_2,c)$\}  & None \\ 

            & 2,2,1,1,1-\{$(c,d),e$\}             &
2,1,1,1,1,1-\{$(a_1,b),(b,c)$\} &  \\

            & 2,2,1,1,1-\{$(b_1,c),(b_2,d)$\}     &
2,1,1,1,1,1-\{$(b,c),(d,e)$\} &  \\ \cline{3-3}

            & 2,2,1,1,1-\{$(a_1,c),(b_1,c)$\}     & 1,1,1,1,1,1-2e &  \\
\cline{2-2}

            & 3,1,1,1,1-\{$(a,b),e$\}             &  &  \\

            & 3,1,1,1,1-\{$(b,c),(d,e)$\}         &  &  \\ \cline{2-2}

            & 4,1,1,1,1-\{$(a_1,b),(a_2,c)$\}     &  &  \\

            & 4,1,1,1,1-\{$(a_1,b),(a_2,b)$\}     &  &  \\  \cline{2-2}

            & 3,2,1,1,1-\{$(a_1,b_1),(a_2,b_2)$\} &  &  \\ \cline{2-2}

            & 2,1,1,1,1-2e                        &  &  \\ \hline

%            &  &  &  \\ 

\end{tabular}

\end{center}

\caption{Intrinsic linking of 2-deficient graphs (cont.)}\label{il2def2}

\end{table}

\begin{table}[ht]

\begin{center}

\begin{tabular}{|l||c|c|c|} \hline

k           & 1    & 2      & 3                                 \\
\hline

knotted     & 8-2e & 5,5-2e & 3,3,3-\{($a_1,b_1$),($a_1,b_2$)\} \\
% & 3,2,2,2-\{$(b,c),e$\}&  3,2,1,1,1-\{($b_1,c$),($b_1,d$)\} \\ 

            &      &        & 3,3,3-\{($a_1,b_1$),($b_2,c_1$)\} \\ \cline{4-4}
% & 3,2,2,2-\{$(a_1,b_1),(a_1,b_2)$\} &  3,2,1,1,1-\{($b_1,c$),($b_2,c$)\} \\ \cline{4-4} \cline{6-6}

            &      &        & 4,4,1-\{($a_1,c$),($b_1,c$)\} \\
%  & 3,2,2,2-\{$(a_1,b_1),(a_2,b_1)$\} &  4,2,1,1,1-\{$(c,d),e$\}  \\

            &      &        & 4,4,1-\{($a_1,b_1$),($b_1,c$)\}   \\ \cline{4-4}
% & 3,2,2,2-\{$(a_1,b_1),(a_2,c_1)$\} & 4,2,1,1,1-\{$(b,c),e$\}  \\ \cline{4-5}

            &      &        & 4,3,2-\{($a_1,b_1$),($a_1,b_2$)\} \\
% & 4,2,2,1-\{$(b,c),e$\} & 4,2,1,1,1-\{($a_1,b_1$),($a_1,c$)\} \\

            &      &        & 4,3,2-\{($a_1,c_1$),($b_1,c_2$)\} \\
% & 4,2,2,1-\{$(c,d),e$\} & 4,2,1,1,1-\{($a_1,b_1$),($a_1,b_2$)\} \\

            &      &        & 4,3,2-\{($a_1,c_1$),($b_1,c_1$)\} \\
% & 4,2,2,1-\{($a_1,b_1$),($a_1,b_2$)\} & 4,2,1,1,1-\{($a_1,c$),($a_1,d$)\} \\ \cline{6-6}

            &      &        & 4,3,2-\{($b_1,c_1$),($b_2,c_1$)\} \\
% & 4,2,2,1-\{($a_1,b_1$),($a_1,d$)\} & 3,3,1,1,1-\{$(b,c),e$\}\\  \cline{5-5}

            &      &        & 4,3,2-\{($a_1,c_1$),($a_1,c_2$)\} \\ 
%& 3,3,2,1-\{$(c,d),e$\} & 3,3,1,1,1-\{$(c,d),e$\}\\

            &      &        & 4,3,2-\{($a_1,b_1$),($b_2,c_1$)\} \\
% & 3,3,2,1-\{$(a,d),e$\} & 3,3,1,1,1-\{$(a_1,b_1),(a_2,b_1)$\}\\ \cline{6-6}

            &      &        & 4,3,2-\{($b_1,c_1$),($b_1,c_2$)\} \\ \cline{4-4}
% & 3,3,2,1-\{($a_1,b_1$),($a_1,b_2$)\} & 3,2,2,1,1-\{$(a,d),e$\}\\ \cline{4-4}

            &      &        & 5,4,1-2e                          \\
% & 3,3,2,1-\{($b_1,c_1$),($b_1,c_2$)\} & 3,2,2,1,1-\{$(c,d),e$\}\\

            &      &        & 5,3,2-2e                          \\
% & 3,3,2,1-\{($a_1,c_1$),($b_1,c_1$)\} & 3,2,2,1,1-\{$(b,c),e$\}\\ 

            &      &        & 4,3,3-2e                          \\
% & 3,3,2,1-\{($a_1,b_1$),($a_2,c_1$)\} & 3,2,2,1,1-\{$(d,e),e$\}\\ 

            &      &        & 4,4,2-2e              \\
% & 3,3,2,1-\{($a_1,c_1$),($a_2,c_1$)\} &3,2,2,1,1-\{$(a_1,b_1),(a_1,b_2)$\}\\ 

\hline

%            &      &        &                                   &  & \\ 

not knotted & 7-2e & n,4-2e & 3,3,3-\{($a_1,b_1$),($b_1,c_1$)\} \\ 
% & 3,2,2,2-\{($a_1,b_1$),($a_1,c_1$)\} & 3,2,1,1,1-\{$(a,b),e$\}\\

            &      &        & 3,3,3-\{($a_1,b_1$),($a_2,b_2$)\} \\ \cline{4-4}
%& 3,2,2,2-\{($a_1,b_1$),($a_2,b_2$)\} &3,2,1,1,1-\{$(a,c),e$\}\\ \cline{4-5}

            &      &        & 4,4,1-\{($a_1,b_1$),($a_1,b_2$)\} \\
% & 4,2,2,1-\{($a_1,b_1$),($a_2,b_1$)\} &3,2,1,1,1-\{$(c,d),e$\}\\

            &      &        & 4,4,1-\{($a_1,b_1$),($a_2,b_2$)\} \\
% & 4,2,2,1-\{($a_1,d$),($a_2,d$)\} & 3,2,1,1,1-\{$(b_1,c),(b_2,d)$\}\\ \cline{6-6}

            &      &        & 4,4,1-\{($a_1,c$),($a_2,c$)\}     \\
% & 4,2,2,1-\{($a_1,b_1$),($a_2,b_2$)\} &4,2,1,1,1-\{($a_1,b_1$),($a_2,b_2$)\}\\

            &      &        & 4,4,1-\{($a_1,b_1$),($a_2,c$)\}   \\ \cline{4-4}
% & 4,2,2,1-\{($a_1,b_1$),($a_2,c_1$)\}&4,2,1,1,1-\{($a_1,b_1$),($a_2,c$)\}\\ \cline{4-4}

            &      &        & 4,3,2-\{($a_1,b_1$),($a_2,b_1$)\} \\
% & 4,2,2,1-\{($a_1,b_1$),($a_2,d$)\} &4,2,1,1,1-\{($a_1,b_1$),($a_2,b_1$)\}\\

            &      &        & 4,3,2-\{($a_1,b_1$),($a_2,b_2$)\} \\
% & 4,2,2,1-\{($a_1,b_1$),($a_1,c_1$)\} &4,2,1,1,1-\{($a_1,c$),($a_2,c$)\}\\ \cline{5-5}

            &      &        & 4,3,2-\{($a_1,c_1$),($a_2,c_1$)\} \\
% & 3,3,2,1-\{($a_1,b_1$),($b_1,c_1$)\} &4,2,1,1,1-\{($a_1,c$),($a_2,d$)\}\\ \cline{6-6}

            &      &        & 4,3,2-\{($a_1,b_1$),($b_1,c_1$)\} \\
% & 3,3,2,1-\{($a_1,b_1$),($a_2,b_2$)\} &3,3,1,1,1-\{($a_1,b_1$),($a_2,b_2$)\}\\ \cline{6-6}

            &      &        & 4,3,2-\{($a_1,b_1$),($a_2,c_1$)\} \\
% & 3,3,2,1-\{($b_1,c_1$),($b_2,c_2$)\} & 3,2,2,1,1-\{($a_1,b_1$),($a_1,c_1$)\}\\ \cline{5-5}

            &      &        & 4,3,2-\{($b_1,c_1$),($b_2,c_2$)\} \\
% & 4,3,1,1-\{($a_1,b_1$),($a_2,b_1$)\} &3,2,2,1,1-\{($a_1,b_1$),($a_2,b_2$)\} \\ \cline{6-6}

            &      &        & 4,3,2-\{($a_1,c_1$),($a_2,c_2$)\} \\
% & 4,3,1,1-\{($a_1,b_1$),($a_2,b_2$)\} &2,2,2,1,1-2e \\

            &      &        & 4,3,2-\{($a_1,b_1$),($a_1,c_1$)\} \\ \cline{4-4}
% & 4,3,1,1-\{($a_1,b_1$),($a_2,c$)\} &n,1,1,1,1-2e\\ \cline{4-4}

            &      &        & 3,3,2-2e                         \\
% & 4,3,1,1-\{($a_1,c$),($a_2,d$)\} &\\

            &      &        & n,2,2-2e                          \\
% & 4,3,1,1-\{($a_1,c$),($a_2,c$)\} &\\ \cline{5-5}

            &      &        & n,3,1-2e                          \\ \hline

%            &      &        &                                   &
% &\\ 

\end{tabular}

\end{center}

\caption{Intrinsic knotting of 2-deficient graphs}\label{ik2def1}

\end{table}

\begin{table}[ht]

\begin{center}

\begin{tabular}{|l||c|c|c|} \hline

k       & \multicolumn{2}{c|}{4} & 5       \\ \hline

knotted & 3,2,2,2-\{$(b,c),e$\}  & 3,3,2,2-2e & 3,2,1,1,1-\{($b_1,c$),($b_1,d$)\} \\ 

        & 3,2,2,2-\{$(a_1,b_1),(a_1,b_2)$\} & 5,2,2,1-2e & 3,2,1,1,1-\{($b_1,c$),($b_2,c$)\} \\ 
\cline{4-4}

        & 3,2,2,2-\{$(a_1,b_1),(a_2,b_1)$\} & 4,3,2,1-2e & 4,2,1,1,1-\{$(c,d),e$\}  \\

        & 3,2,2,2-\{$(a_1,b_1),(a_2,c_1)$\} & 3,3,3,1-2e & 4,2,1,1,1-\{$(b,c),e$\}  \\ \cline{2-2}

        & 4,2,2,1-\{$(b,c),e$\}             & 5,3,1,1-2e & 4,2,1,1,1-\{($a_1,b_1$),($a_1,c$)\} \\

        & 4,2,2,1-\{$(c,d),e$\}             & 4,4,1,1-2e & 4,2,1,1,1-\{($a_1,b_1$),($a_1,b_2$)\} \\
 
        & 4,2,2,1-\{($a_1,b_1$),($a_1,b_2$)\} & & 4,2,1,1,1-\{($a_1,c$),($a_1,d$)\} \\ \cline{4-4}

        & 4,2,2,1-\{($a_1,b_1$),($a_1,d$)\} &  & 3,3,1,1,1-\{$(b,c),e$\}\\  \cline{2-2}

        & 3,3,2,1-\{$(c,d),e$\} & & 3,3,1,1,1-\{$(c,d),e$\}\\

        & 3,3,2,1-\{$(a,d),e$\} & & 3,3,1,1,1-\{$(a_1,b_1),(a_2,b_1)$\}\\ \cline{4-4}

        & 3,3,2,1-\{($a_1,b_1$),($a_1,b_2$)\} & & 3,2,2,1,1-\{$(a,d),e$\}\\

        & 3,3,2,1-\{($b_1,c_1$),($b_1,c_2$)\} & & 3,2,2,1,1-\{$(c,d),e$\}\\

        & 3,3,2,1-\{($a_1,c_1$),($b_1,c_1$)\} & & 3,2,2,1,1-\{$(b,c),e$\}\\ 

        & 3,3,2,1-\{($a_1,b_1$),($a_2,c_1$)\} & & 3,2,2,1,1-\{$(d,e),e$\}\\ 

        & 3,3,2,1-\{($a_1,c_1$),($a_2,c_1$)\} & & 3,2,2,1,1-\{$(a_1,b_1),(a_1,b_2)$\}\\ 

        & 3,3,2,1-\{($a_1,c_1$),($b_1,c_2$)\} & & 3,2,2,1,1-\{$(a_1,b_1),(a_2,b_1)$\}\\  \cline{2-2}

        & 4,3,1,1-\{$(b,c),e$\}               & & 3,2,2,1,1-\{$(a_1,b_1),(a_2,c_1)$\}\\ \cline{4-4}

        & 4,3,1,1-\{$(c,d),e$\}               & & 2,2,2,2,1-2e\\ 

        & 4,3,1,1-\{($a_1,b_1$),($a_1,b_2$)\} & & 5,2,1,1,1-2e\\ 

        & 4,3,1,1-\{($a_1,b_1$),($a_1,c$)\}   & & 4,3,1,1,1-2e\\ 

        & 4,3,1,1-\{($a_1,c$),($a_1,d$)\}     & & 4,2,2,1,1-2e\\ \cline{2-2}

        & 4,2,2,2-2e                          & & 3,3,2,1,1-2e \\ \hline

%            &      &        &                                   &  & \\ 

not         & 3,2,2,2-\{($a_1,b_1$),($a_1,c_1$)\} & 3,3,1,1-2e & 3,2,1,1,1-\{$(a,b),e$\}\\

knotted     & 3,2,2,2-\{($a_1,b_1$),($a_2,b_2$)\} & 2,2,2,2-2e & 3,2,1,1,1-\{$(a,c),e$\}\\
\cline{2-2}

            & 4,2,2,1-\{($a_1,b_1$),($a_2,b_1$)\} & 3,2,2,1-2e & 3,2,1,1,1-\{$(c,d),e$\}\\

            & 4,2,2,1-\{($a_1,d$),($a_2,d$)\}     & n,2,1,1-2e & 3,2,1,1,1-\{$(b_1,c),(b_2,d)$\}\\
\cline{4-4}

            & 4,2,2,1-\{($a_1,b_1$),($a_2,b_2$)\} & & 4,2,1,1,1-\{($a_1,b_1$),($a_2,b_2$)\}\\

            & 4,2,2,1-\{($a_1,b_1$),($a_2,c_1$)\} & & 4,2,1,1,1-\{($a_1,b_1$),($a_2,c$)\}\\

            & 4,2,2,1-\{($a_1,b_1$),($a_2,d$)\}   & & 4,2,1,1,1-\{($a_1,b_1$),($a_2,b_1$)\}\\

            & 4,2,2,1-\{($a_1,b_1$),($a_1,c_1$)\} & & 4,2,1,1,1-\{($a_1,c$),($a_2,c$)\}\\ \cline{2-2}

            & 3,3,2,1-\{($a_1,b_1$),($b_1,c_1$)\} & & 4,2,1,1,1-\{($a_1,c$),($a_2,d$)\}\\ \cline{4-4}

            & 3,3,2,1-\{($a_1,b_1$),($a_2,b_2$)\} & & 3,3,1,1,1-\{($a_1,b_1$),($a_2,b_2$)\}\\
\cline{4-4}

            & 3,3,2,1-\{($b_1,c_1$),($b_2,c_2$)\} & & 3,2,2,1,1-\{($a_1,b_1$),($a_1,c_1$)\}\\
\cline{2-2}

            & 4,3,1,1-\{($a_1,b_1$),($a_2,b_1$)\} & & 3,2,2,1,1-\{($a_1,b_1$),($a_2,b_2$)\} \\
\cline{4-4}

            & 4,3,1,1-\{($a_1,b_1$),($a_2,b_2$)\} & & 2,2,2,1,1-2e \\

            & 4,3,1,1-\{($a_1,b_1$),($a_2,c$)\}   & & n,1,1,1,1-2e\\ 

            & 4,3,1,1-\{($a_1,c$),($a_2,d$)\} & & \\

            & 4,3,1,1-\{($a_1,c$),($a_2,c$)\} & & \\ \hline

\end{tabular}

\end{center}

\caption{Intrinsic knotting of 2-deficient graphs (cont.)}\label{ik2def2}

\end{table}

Next, we'll argue that the ``not linked" graphs are in fact not intrinsically linked.

$K_6$ (or $K_{1,1,1,1,1,1}$) and $K_{3,3,1}$ are Petersen graphs and
therefore minor minimal with respect to intrinsic linking
\cite{RST}, so $K_6$-e and $K_{3,3,1}$-e are not intrinsically
linked by definition of minor minimal.

$K_{n,3}$-e, $K_{3,2,2}$-e, $K_{n,2,1}$-e, $K_{2,2,2,1}$-e,
$K_{n,1,1,1}$-e, and $K_{2,1,1,1,1}$-e are all not intrinsically linked
before the edge is removed; so, naturally, we can remove an edge and still
have a non-intrinsically linked graph.

$K_{3,2,1,1}$-(c,d) is equivalent to $K_{3,2,2}$ which is not
intrinsically linked.

The remaining graphs all have a vertex that is connected to every other
vertex, and the removal of that vertex results in a planar graph. So, by
Lemma \ref{sachs}, none are intrinsically linked.
\qed

\subsection{Intrinsic knotting}

\begin{thm}

The 1-deficient graphs are classified with respect to intrinsic  knotting
according to Table \ref{ik1def}.
\end{thm}

\Pf

%\noindent%
%\underline{Knotted:}

Let us demonstrate that the graphs labeled ``knotted" in Table
\ref{ik1def} are in fact intrinsically knotted. 

By Theorem~\ref{thmIK}, 
$K_8$-e, $K_{2,2,2,2,1}$-e, $K_{3,2,1,1,1,1}$-e,
$K_{2,2,2,1,1,1}$-e, and $K_{2,1,1,1,1,1,1}$-e are intrinsically knotted.
Any $1$-deficient graph of eight or more parts will contain
$K_8$-e as a minor and therefore be intrinsically knotted.

$K_{5,5}$-e is intrinsically knotted as proved in \cite{Sh}. 

$K_{3,3,3}$-e and $K_{4,4,1}$-e have the graph $H_9$ of \cite{KS} (derived from
$K_7$ by two triangle-Y moves) as 
a minor and are therefore intrinsically knotted. By Lemma~\ref{combparts},
$K_{3,3,2,1}$-e, $K_{3,3,1,1,1}$-e,
and $K_{3,2,2,1,1}$-e all have
$K_{3,3,3}$-e as a minor, so they too are intrinsically knotted.

$K_{4,3,2}$-e is intrinsically knotted. Depending on which edge is removed,
the graph will either have $H_9$ of \cite{KS} or $K_{3,3,1,1}$ as a minor.
Consequently, by Lemma~\ref{combparts},
$K_{3,2,2,2}$-e, $K_{4,2,2,1}$-e,
$K_{4,3,1,1}$-e, $K_{4,2,1,1,1}$-e, and
$K_{4,1,1,1,1,1}$-e are intrinsically knotted.

Finally, $K_{3,2,1,1,1}$-(b,c), $K_{2,2,1,1,1,1}$-(b,c), and
$K_{3,1,1,1,1,1}$-(b,c), all have $K_{3,3,1,1}$ as a minor and are therefore
intrinsically knotted.

Now we'll show that the ``not knotted" graphs are not in fact
intrinsically knotted.

$K_7$ was shown to be minor minimally intrinsically knotted in
\cite{CG} and $K_{3,3,1,1}$ was shown to be minor minimally
intrinsically knotted in \cite{F1}, so $K_7$-e (or
$K_{1,1,1,1,1,1,1}$-e) and $K_{3,3,1,1}$-e are not intrinsically knotted.

In \cite{Fl}, the authors show that $K_{4,4}$ is not
intrinsically knotted.  In fact, the proof really shows that $K_{n,4}$ is not
intrinsically knotted, so, naturally, $K_{n,4}$-e is not intrinsically
knotted either. Similarly, $K_{3,3,2}$-e, $K_{n,2,2}$-e, $K_{n,3,1}$-e,
$K_{2,2,2,2}$-e, $K_{3,2,2,1}$-e, $K_{n,2,1,1}$-e, $K_{2,2,2,1,1}$-e,
$K_{n,1,1,1,1}$-e, and $K_{2,1,1,1,1,1}$-e are not intrinsically knotted.

$K_{3,2,1,1,1}$-(c,d) is equivalent to $K_{3,2,2,1}$.

The remaining 5 graphs each have 2 vertices connected to every other
vertex. If those 2
vertices are deleted, the resulting graph is planar, so by Lemma \ref{flemming},
they are not intrinsically knotted.
\qed

\begin{table}[ht]

\begin{center}

\begin{tabular}{|l||c|c|c|} \hline

k             & 6                                   & 7 & $\geq$8 \\
\hline

knotted       & 2,2,1,1,1,1-\{($b_1,c$),($b_2,c$)\} &
2,1,1,1,1,1,1-\{($a_1,b$),($a_1,c$)\} & all \\ 

              & 2,2,1,1,1,1-\{($b_1,c$),($b_1,d$)\} &
2,1,1,1,1,1,1-\{($a_1,b$),($a_2,b$)\} &  \\ 

              & 2,2,1,1,1,1-\{($a_1,d$),($b_1,c$)\}     &
2,1,1,1,1,1,1-\{($a_1,b$),($c,d$)\} &  \\ \cline{2-2}

              & 3,1,1,1,1,1-\{($b,c$),($c,d$)\}           &
2,1,1,1,1,1,1-\{($b,c$),($c,d$)\} &  \\ \cline{2-3}

              & 4,1,1,1,1,1-\{$(b,c),e$\} & 3,1,1,1,1,1,1-2e &  \\ 

              & 4,1,1,1,1,1-\{($a_1,b$),($a_1,c$)\}      &
2,2,1,1,1,1,1-2e &  \\ \cline{2-2}

              & 3,2,1,1,1,1-\{$(a,c),e$\}   &  &  \\ 

              &  3,2,1,1,1,1-\{$(b,c),e$\}  &  &  \\ 

              & 3,2,1,1,1,1-\{$(c,d),e$\}    &  &  \\ 

              & 3,2,1,1,1,1-\{$(a_1,b_1),(a_1,b_2)$\}    &  &  \\ 

              & 3,2,1,1,1,1-\{$(a_1,b_1),(a_2,b_1)$\} &  &  \\ \cline{2-2}

              & 2,2,2,1,1,1-2e &  &  \\ 

              & 5,1,1,1,1,1-2e &  &  \\ 

              & 4,2,1,1,1,1-2e &  &  \\ 

              & 3,3,1,1,1,1-2e &  &  \\ \hline

not           & 2,2,1,1,1,1-\{$(a,b),e$\}     &
2,1,1,1,1,1,1-\{($a_1,b$),($b,c$)\} & none \\ 

knotted       &  2,2,1,1,1,1-\{$(c,d),e$\}     &
2,1,1,1,1,1,1-\{($b,c$),($d,e$)\} &  \\ 

              & 2,2,1,1,1,1-\{$(b_1,c),(b_2,d)$\} &
2,1,1,1,1,1,1-\{($a_1,b$),($a_2,c$)\} &  \\ \cline{3-3}

              & 2,2,1,1,1,1-\{$(a_1,c),(b_1,c)$\} & 1,1,1,1,1,1,1-2e & 
\\  \cline{2-2}

              & 3,1,1,1,1,1-\{$(a,b),e$\} &  &  \\ 

              & 3,1,1,1,1,1-\{$(b,c),(d,e)$\} &  &  \\ \cline{2-2}

              & 4,1,1,1,1,1-\{($a_1,b$),($a_2,b$)\} &  &  \\ 

              & 4,1,1,1,1,1-\{($a_1,b$),($a_2,c$)\} &  &  \\ \cline{2-2}

              & 3,2,1,1,1,1-\{($a_1,b_1$),($a_2,b_2$)\} &  &  \\ \cline{2-2}

              & 2,1,1,1,1,1-2e &  &  \\  \hline

\end{tabular}

\end{center}

\caption{Intrinsic knotting of 2-deficient graphs (cont.)}\label{ik2def3}

\end{table}

\section{$2$-deficient graphs}

In this section we give the classification of 
$2$-deficient graphs with respect to 
intrinsic linking and knotting in Tables~\ref{il2def1}, \ref{il2def2},
\ref{ik2def1}, \ref{ik2def2}, and \ref{ik2def3}. The argument justifying this classification
is very similar to that presented in the last section and we will not 
include it here. See~\cite{MOR} for a full account.

\Nt We will expand the notation that we created in the last section
by adding subscripts to the vertices. For example,
if we are removing two edges  between parts A and B of $K_{4,3,1}$,
there are 3 cases to be considered: deleting two edges that share a
vertex from part A ($K_{4,3,1}$-\{$(a_1,b_1),(a_1,b_2)$\}), deleting two
edges that share a vertex from part B
($K_{4,3,1}$-\{$(a_1,b_1),(a_2,b_1)$\}), and deleting two edges that
share no vertices ($K_{4,3,1}$-\{$(a_1,b_1),(a_2,b_2)$\}). Also, in some
cases, if we delete a particular edge we can then delete any other 
without affecting the classification of the resulting graph. For
example, if we delete edge (b,c) from $K_{4,3,1}$, we can then delete any
other and still have an intrinsically linked graph; we will denote graphs
obtained in this way by
$K_{4,3,1}$-\{$(b,c),e$\}. Furthermore, in some cases, we can delete any
2 edges and still have an intrinsically linked graph;
for example, $K_7$ will be
intrinsically linked no matter what 2 edges we delete. We denote these
graphs $K_7$-2e.

\section{Graphs on $8$ vertices}
In this section we provide a classification of intrinsically knotted
graphs on $8$ vertices. This classification was recently independently presented 
by Blain, Bowlin, Fleming, Foisy, Hendricks, and LaCombe~\cite{Fl}.

Graphs on $8$ vertices are subgraphs of $K_8$. We will examine in turn
subgraphs which are obtained by removing $1$, $2$, $3$, $\ldots$ edges 
from $K_8$.  We have already noted that
$K_8$, $K_8 - e$, and both $K_8-2e$ graphs are intrinsically
knotted. Of the five graphs $K_8-3e$, the three which are intrinsically 
knotted can all be obtained by removing two edges from
$K_{2,1,1,1,1,1,1}$. (See classification of $2$-deficient graphs
above. Note that $2,1,1,1,1,1,1-\{(a,b),(c,d)\}$ and
$2,1,1,1,1,1,1 - \{(b,c),(c,d) \}$ are the same graph.)

There are $11$ graphs of the form $K_8-4e$. Seven of these are not
knotted as they can be realized by removing an edge from one of the two
unknotted
$K_8-3e$ graphs. The remaining four are intrinsically knotted. Three of
these four are  of the form $K_{2,2,1,1,1,1} - 2e$. The fourth is
obtained by removing
$4$ edges all incident to the same vertex. This graph is intrinsically
knotted as it contains $K_7$ as a minor.

There are $24$ graphs $K_8-5e$. All but $4$ of these are  not
knotted as they are minors of an unknotted $K_8-4e$. The four
intrinsically knotted $K_8-5e$'s are perhaps most easily described in
terms of their complementary graphs.
\begin{figure}[ht]
\begin{center}
\includegraphics[scale=0.2]{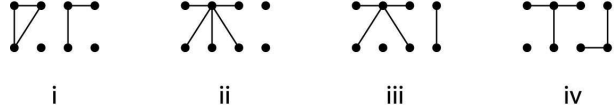}
%\centerline{\includegraphics{fig9.eps}}
%\vspace{0.7 in}
%\hspace{-3 in}
%\special{illustration fig9.eps scaled 0.5}

\caption{Intrinsically knotted $K_8-5e$'s}\label{fig5def}
\end{center}
\end{figure} 
In Figure~\ref{fig5def}, Graph i is
$K_{3,2,1,1,1} - (b,c)$, or, equivalently, $K_{3,1,1,1,1,1} - \{(b,c),(c,d)\}$.
The other three graphs in the figure have $K_7$ as a minor and are,
therefore, intrinsically knotted.

Of the $56$ $K_8 - 6e$ graphs, all but $6$ are minors of unknotted
$K_8 - 5e$'s. 
\begin{figure}[ht]
\begin{center}
\includegraphics[scale=0.1]{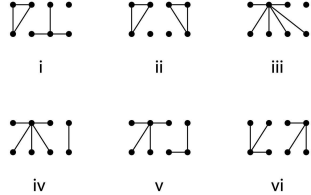}
%\centerline{\includegraphics{fig10.eps}}
%\vspace{1.4 in}
%\hspace{-2 in}
%\special{illustration fig10.eps scaled 0.5}

\caption{Intrinsically knotted $K_8-6e$'s}\label{fig6def}
\end{center}
\end{figure} 
In Figure~\ref{fig6def}, Graph i is $K_{3,2,1,1,1} -
\{(b_1,c),(b_1,d)\}$ while Graph ii is $K_{3,3,1,1}$, or, equivalently,
$K_{3,2,1,1,1} - \{(b_1,c),(b_2,c)\}$. The remaining four graphs are
obtained by splitting a vertex of $K_7$  and are therefore intrinsically
knotted.

Only $2$ of the $K_8-7e$ graphs are not minors of some
unknotted $K_8 -6e$. One of these two is $H_8$ \cite{KS}, the graph obtained by a single
triangle-Y exchange on $K_7$. The other is
$K_7$ with one additional vertex. These are both intrinsically knotted. 
Moreover, $K_7$ and $H_8$ are minor minimal \cite{KS}. Thus any subgraph of $K_8$
obtained by removing $8$ or more edges is not intrinsically knotted. 

In total then, there are twenty intrinsically knotted graphs on $8$
vertices. 

%\section*{Acknowledgements}
%This work was completed in part during REUs at CSU, Chico during the
%summers of 2003 through 2005 supported by NSF REU award \#0354174 and
%the MAA's Strengthening Underrepresented Minority Mathematics Achievement
%Program, with funding from the NSA, the
%NSF, and the Moody's Foundation. Additional funding came from CSU,
%Chico Research Foundation.

\end{document}